\documentclass[12pt]{amsart}
\usepackage{a4wide}
\usepackage{amssymb}
\usepackage{wrapfig,epsfig}
\usepackage{color}
\usepackage{mathrsfs}
\usepackage{pdfsync}

\usepackage{hyperref}

\numberwithin{equation}{section}

\newcommand\dela[1]{}
\newcommand\delb[1]{#1}

\newcommand\sol{\mathrm{sol}}
\newcommand{\embed}{\hookrightarrow}
\newcommand{\toup}{\nearrow}

\theoremstyle{definition}
\newtheorem{definition}{Definition}[section]

\newtheorem{remark}[definition]{Remark}

\theoremstyle{plain}

\newtheorem{theorem}[definition]{Theorem}
\def\REQ#1{{\rm (\ref{#1})}}

\def\d{{\rm d}}
\def\e{{\rm e}}

\begin{document}

\title{Ergodicity for stochastic equations of  Navier--Stokes type}
\author{Zdzis{\l}aw Brze{\'z}niak}
\address{Department of Mathematics, The University of York, Heslington, York YO10 5DD, UK.}
\email{zb500@york.ac.uk}
\author{Tomasz Komorowski}
\address{Institute of Mathematics, Polish Academy of Sciences, \'Sniadeckich, 8, 00-956, Warsaw, Poland.}
\email{tkomorowski@impan.pl}
\author{Szymon Peszat}
\address{Faculty of Mathematics and Computer Science,  Jagiellonian University in Krakow t{\L}ojasiewicza, 30-348 Krak{\'o}w, Poland } 

\email{szymon.peszat@im.uj.edu.pl}



\begin{abstract} In the first part of the note we analyze the long time behaviour of  a  two dimensional stochastic Navier--Stokes equations system on a torus  with a degenerate, one dimensional noise.  In particular, for some initial data and noises we identify the invariant probability measure for the system and give a sufficient condition under which it is unique and stochastically stable. In the second part of the note,
we consider a simple example of a finite dimensional system of stochastic differential equations  driven by a one dimensional Wiener process with a drift, that displays some similarity with the stochastic N.S.E., and investigate its ergodic properties depending on the strength of the drift.  If the latter is sufficiently small and lies below a critical threshold,
then the system admits a unique invariant probability measure which is Gaussian.  If, on the other hand, the  strength of the noise drift  is larger than the  threshold, then in addition to a Gaussian invariant probability measure, there exist another one. In particular, the generator of
the system is not hypoelliptic.
\end{abstract}

\maketitle

\section{Introduction}

Study of ergodic properties of dynamical (inclusive random) systems is of profound importance from both applied and theoretical standpoints. Two examples of such properties are the existence and uniqueness (or possibly non-uniqueness) of  invariant probability measures. These are often linked to the not-yet fully  explained aspects of turbulence such as e.g.  the rigorous proof of the form of the Kolmogorov spectrum.

In the case of stochastic hydrodynamics, i.e. for the stochastic Navier--Stokes equations of the following form
\begin{equation}\label{eqn_NSEs}
\partial_t u + (u\cdot \nabla)u = \left(\mu\Delta u -\nabla p +f \right) +\xi, \qquad \nabla \cdot u=0,
\end{equation}
where $u$ is the velocity field, $p$ is the pressure scalar (both unknown), $f$ is the external force acting on tyhe fluid and $\xi$ is a noise,
the first results in those directions are due to Flandoli \cite{Flandoli},
 showed the existence of an invariant probability measure for the 2D Navier--Stokes  equations (NSE) driven by an additive Gaussian  noise. The question of the uniqueness of an invariant probability measure   for the same system  has been later addressed by   Flandoli and Maslowski \cite{Flandoli+Maslowski_1995}, Ferario \cite{Ferrario_1997} and  E, Mattingly and Sinai \cite{EMS}. The first two papers assumed that noise $\xi$ has been sufficiently non-degenerate (which had to be counterbalanced by a requirement that it is sufficiently spatially regular to ensure the solvability of the system \eqref{eqn_NSEs}).  The direction of research initiated in first two of these papers   has been followed by Da Prato and Debussche in \cite{DaP+Deb_2002} who proved the unique solvability of \eqref{eqn_NSEs} (in Besov spaces of negative order) and the existence of an invariant probability measure for the stochastic Stokes equation (on a 2D torus) when $\xi$ is the space-time noise; see also  the paper by Albeverio and Ferrario \cite{Alb+Fer_2002}.

 The paper \cite{EMS} by  E, Mattingly and Sinai  looked at the question of uniqueness of an  invariant probability measure in the case of a degenerate noise, which happens to be mathematically more challenging than the non-degenerate noise. In this case the corresponding Markov process is
 is only a Feller rather than strong Feller. This case was also studied by Mattingly in \cite{Mat99,Mat03} in the case the external force $f$ is equal to zero and in \cite{GMR17,KS18,BKS20} for nonzero force $f$.
 The culminating work on this topic is due to Hairer and Mattingly \cite{HM} who, using a new concept of an asymptotically strong Feller semigroup, proved that the Markov process generated by the   stochastic NSEs  on a 2-d torus has a unique invariant probability measure provided  the Gaussian perturbation is of mean $0$ and acts on at least two  modes that are of different length and whose integer linear combinations generate the two dimensional integer lattice. Such a system can be called a hypoelliptic. Later on Friedlander et al. \cite{F+GH+V_2016} and   Andreis et al. \cite{A+B+C+F+P_2016}, proved that the hypoellipticity still holds for certain  stochastic inviscid dyadic models and hence such models have a unique invariant probability measure even if the centered noise acts only on a single mode.

  It is still an open question whether similar properties hold in the presence of a large deterministic force, i.e. when the noise in not centered and its mean is large. For instance the method from \cite{HM}  still works when the force is small so that the corresponding deterministic system has a unique stationary solution which is exponentially stable.  Another open question is whether whether similar properties hold  when the noise is more degenerate than the noise considered in the paper \cite{HM}.

  The modest aim of this note is to prove that for a certain finite dimensional system modelling the true SNSE, introduced by Minea in  \cite{Minea},  such a result is not true. To be precise, in Theorem \ref{T34} we show that if  $\kappa  >\lambda_1\min\{\lambda_2, \lambda_3\}$,  then  the stochastic system \eqref{E36}, i.e. \eqref{E36b}, has at least two  invariant probability measures, and,  since any  convex combination of these invariant probability measures is also an invariant probability measure,  the stochastic system \eqref{E36}, i.e.  \eqref{E36b}, has infinitely many invariant probability measures.

One of the measures, denoted by $\nu_{\sigma,\kappa}$,  is  Gaussian.  This measure  is also  the unique  invariant probability measure for the corresponding stochastic "Stokes  system" \eqref{E36d}. Let us finish this paragraph by recalling that the set of stationary solutions for the corresponding deterministic system \eqref{E39b}  has quite a  complicated structure. Thus the  present note  shows that this also could be  the case for its stochastic perturbation.

The paper is organized as follows. In Section \ref{sec2} we recall the basic functional analytic setting used for the evolution equations considered the paper. The main results are formulated in Section \ref{sec3}. Their proofs are presented in Section \ref{sec4}.

\section{Preliminaries}
\label{sec2}

Let  $C^{\infty}_{\sol}$ be the set of all  infinitely differentiable $\mathbb{R}^2$-valued mappings $v=(v_1,v_2)$ defined on a  two dimensional standard torus $ \mathbb{T}^2$  such that ${\rm div}\, v=0$ and $\int_{\mathbb{T}^2} v(x)\d x=0$.  Given $q\ge 1$ let us denote by  $L^q_\sol$ the closure of $C^{\infty}_{\sol}$  in $L^q(\mathbb{T}^2):=L^q({\mathbb{T}^2}, \mathcal{B}(\mathbb{T}^2), \d x;{\mathbb{R}}^2)$. We denote  the space $L^2_\sol$ by $H$ and by $\pi$  we  denote the orthogonal projection $\pi \colon L^2 \to H$.
The scalar products and the norms on spaces $L^2 $ and $H$  are denoted, respectively,  by  $(\cdot ,\cdot )$ and $|\cdot |$.

Let $A$ be the Stokes operator on $H$, which is the self-adjoint  operator obtained by the closure  of $- \pi \Delta$ considered on $C^{\infty}_{\sol}$.  Note that $A$ is strictly positive $A^{-1}$ exists and is  a compact symmetric operator.
We will denote by $\lambda_0$ the smallest eigenvalue of the operator $A$.
For $r\ge 0$ we define $H^r_\sol$ as the domain of  $A^{r/2}$.  The space $H^r_\sol$ is equipped with the graph norm  $|\cdot |_r:=|A^{r/2}\cdot |$.  Clearly, $H=H^0_\sol$ and $|\cdot |=|\cdot |_0$.  It is customary to denote  the space $H^1_\sol$ by $V$ and the norm on $V$ by $\|\cdot \|$. Note that the following Poincar{\`e} type inequality holds
\[
\| u  \| \geq \lambda_0^{1/2} |u|, \;\; u \in V.
\]
In particular, $\lambda_0^{1/2}$ is equal to the norm of the canonical embedding $V \embed H$.

Let $V^\prime$ be the dual of $V$; then $H$ can be identified with a subspace of $V^\prime$ and $V\subset H\equiv H^\prime \subset V^\prime$.  Note that by the Sobolev embedding in dimension $d=2$, the   space $V$ is compactly embedded into  $L^q(\mathbb{T}^2)$, for any $q\in [1,+\infty)$.

Define
$$
b(u,v,w)=( (u \cdot \nabla) v,  w)=\sum_{i,j=1}^2 \int_{\mathbb{T}^2} u^i(x) \frac{\partial v^j}{\partial x_{i}}(x) w^j(x)\d x
$$
whenever the integrals make sense.  Below we list some  well-known inequalities  (see \cite[pp. 108]{Temam}) concerning the triple product form $b(\cdot,\cdot,\cdot)$,
\begin{equation}\label{E21}
|b(u,v,w)|\le C\times \left\{
\begin{array}{l}
|u|^{1/2} \|u\|^{1/2}\|v\|^{1/2} |Av|^{1/2}|w|, \quad u\in V, \ v\in H^2_\sol,\ w\in H,\\
|u|^{1/2} |Au|^{1/2}\|v\| |w|, \quad u\in H^2_\sol, \ v\in V,\ w\in H,\\
|u| \|v\|  |Aw|^{1/2}|w|^{1/2}, \quad u\in H, \ v\in V,\ w\in H^2_\sol,\\
|u| ^{1/2}\|u\| ^{1/2} \|v\| |w|^{1/2}\|w\|^{1/2}, \quad u, \, v,\, w\in V,\\
\end{array}
\right.
\end{equation}
where $C>0$ is an appropriate constant.  Hence the bilinear operator $B$ defined  by
$$
(B(u,v),w)=b(u,v,w)
$$
maps continuously $V\times V$ into $V^\prime$.  We will use the facts
\begin{align}\label{E22}
b(u,v,w)&=-b(u,w,v),\qquad u,v,w\in V,\\
\label{E23}
(B(u,v),v)&= 0,\qquad u,v\in V,\\
\label{E24}
B(e,e)&=0 \qquad \text{for  any eigenvector  $e$ of $A$,}\\
\label{E25}
(B(v,v),Av)&=0,\qquad  v\in H^2_\sol.
\end{align}

It is worth noting that \REQ{E24} and \REQ{E25}  hold only in the periodic $2$-dimensional case,  whereas  \REQ{E22} and \REQ{E23} hold  both in dimensions $2$ and $3$ and also in the case of  the Dirichlet boundary conditions.
The proof of property \eqref{E24} can be found in \cite{Flandoli+Gozzi_1998}.

\section{Main results} \label{sec3}

\subsection{The case of Navier--Stokes equations}
Let $e\not =0$ be a normalised  eigenvector of the Stokes operator, and let $\lambda>0$ be its corresponding eigenvalue; that is $Ae=\lambda e$.
Let us observe that $\|e\|^2 =\lambda$.  Let  $\mu>0$ be  the viscosity of the fluid, $\kappa \in \mathbb{R}$,  $\sigma >0$, and let $W$ be a standard real-valued Wiener process.  Given $v\in H$, we denote by $u(\cdot;v)$ the solution of   the following   Navier--Stokes equations on $\mathbb{T}^2$,
\begin{equation}\label{E31}
\d u = \left(\mu A u + B(u,u)\right)\d t+ \kappa e\d t + \sigma e\d W, \qquad u(0;v)=v.
\end{equation}

Let $z(\cdot;v)$ be the solution of  the  corresponding Stokes  equation
\begin{equation}\label{E32}
\d z = \mu Az \d t +  \kappa e\d t +\sigma e\d W, \qquad z(0;v)=v.
\end{equation}
Process $z$ is usually called an Ornstein--Uhlenbeck process.

Let $\nu_{\lambda,\mu,e, \sigma, \kappa }$ be the law of an $V$-valued random variable 
$$
\frac{\kappa}{\lambda \mu} e+ \sigma  \int_0^{+\infty} \e ^{-\lambda \mu t}\d W(t)e,
$$
 i.e.
\begin{equation}\label{E33}
\nu_{\lambda,\mu,e, \sigma, \kappa }:= \mathcal{L}\left(\frac{\kappa}{\lambda \mu} e+ \sigma  \int_0^{+\infty} \e ^{-\lambda \mu t}\d W(t)e\right)= \mathcal{L}\bigl(\eta\, e\bigr),
\end{equation}
where $\eta$ is real gaussian variable with mean $\frac{\kappa}{\lambda \mu}$ and variance $\frac{\sigma^2}{2\lambda \mu}$.
It is known, see e.g. Theorem 11.7 from \cite{DZ1}, that $\nu_{\lambda,\mu,e, \sigma, \kappa }$ is a  unique (and consequently   ergodic)  invariant probability measure to the Langevin equation \eqref{E32}. Let us recall that %
\begin{equation}
\label{E34}
 \int_{V}\|z\|^2 \nu_{\lambda,\mu,e, \sigma, \kappa }(\d z)=\left( \frac{\kappa^2}{\lambda^2\mu^2}+\frac{\sigma^2}{2\lambda \mu}\right)\|e\|^2 = \left( \frac{\kappa^2}{\lambda^2\mu^2}+\frac{\sigma^2}{2\lambda \mu}\right)\lambda.
\end{equation}

\begin{theorem}\label{T31}
\noindent
\begin{itemize}
\item[$(i)$] For any initial data $v\in H$ there exists a sequence $t_n\toup +\infty$ such that the following sequence of  probability measures on $H$, $$
\mathcal{L}\left(\frac{1}{t_n}\int_0^{t_n} u(s;v)\d s\right)
$$
 converges weakly to a certain probability measure on $(H,\mathcal{B}(H))$. Consequently, by the Krylov--Bogolyubov theorem, the stochastic NSE  \eqref{E31} has at least one invariant probability measure. \\
\item[$(ii)$] If the initial data is of the form $v=ae$, where $a\in \mathbb{R}$ and $e$ is a normalised   eigenvector of $A$, then $u(t;v)=z(t;v)$ for all $t\ge 0$. Consequently the invariant  measure $\nu_{\lambda,\mu,e,\sigma,\kappa}$ of  the Langevin equation   \REQ{E32} is an invariant probability measure for the stochastic Navier--Stokes equations \eqref{E31}.\\
\item[$(iii)$]  Let $C$ be the constant appearing in \REQ{E21}. If $\mu$, $\lambda$, $e$, $\sigma$, and  $\kappa$  are such that
\begin{equation}\label{E35}
  C^2 \left( \frac{\kappa^2}{\lambda\mu^4}+\frac{\sigma^2}{2 \mu^3}\right)< \lambda_0,
\end{equation}
then  for any  $v\in H$, the laws  $\mathcal{L}(u(t;v))$ converges weakly, as $t\to +\infty$,  to the invariant probability measure  $\nu_{\lambda,\mu,e,\sigma,\kappa}$ to the Langevin equation \eqref{E32}. Therefore $\nu_{\lambda,\mu,e,\sigma,\kappa}$ is the unique invariant probability measure to  the stochastic Navier--Stokes equations \eqref{E31}.
\end{itemize}
\end{theorem}

\begin{remark}\label{R32}
Condition \REQ{E35}  implies  that both $\kappa$ and $\sigma$ are small enough. In fact, as it can be seen from \REQ{E33},  given $\lambda,\mu,e$, condition    \REQ{E35} is violated when  either   $\kappa$ or   $\sigma$ is  large.
\end{remark}

\begin{remark}\label{R33}
 The existence of an invariant probability measure given in part (i) of Theorem \ref{T31} is classical even for more general stochastic 2D Navier--Stokes equations,  see e.g. Flandoli \cite{Flandoli} for the case of SNSEs in a bounded domain with Dirichlet boundary and \cite{Brz+Li_2006} for   the case of SNSEs  in unbounded Poincar\`e domains.
We present a short proof of this fact  for the sake of completeness.
\end{remark}

\begin{remark}\label{R33a}
The fact that an invariant probability measure for the stochastic Stokes equation (on a 2D torus) \REQ{E31}   driven by a canonical cylindrical Wiener process on $H$ is also  an invariant probability measure for the corresponding stochastic Navier--Stokes equations \REQ{E32} is known, see e.g. the paper \cite{DaP+Deb_2002} by Da Prato and Debussche, where this statement is made rigorous, and  also  the paper by Albeverio and Ferrario \cite{Alb+Fer_2002}.
\end{remark}

\begin{remark}\label{R33b}
    The result given in part (i) of Theorem \ref{T31} is
(the existence of an invariant probability measure for the stochastic NSEs  and the applicability of the Krylov--Bogolyubov) is classical and can be traced back at least to the
paper by  Flandoli \cite[Section 3.3]{Flandoli}.
  The result given in part (iii) of Theorem \ref{T31} is    known in the case  when  external force $f=\kappa e$ is equal to $0$, see e.g. the paper \cite{EMS} by E,  Mattingly and  Sinai. Namely,  \cite[Theorem 1]{EMS} shows uniqueness of the invariant probability measure for
the stochastic NS equation if $\kappa = 0$  and $\sigma^2/\nu^3 < C$. Note that, as explained above,
the result given in part (iii) of Theorem 3.1 is covered completely in later works \cite{GMR17,KS18,BKS20}.
\end{remark}

\subsection{Simplified  Navier--Stokes Equations}\label{subsec-simplified NSEs}

Consider the following $\mathbb R^3$-valued system of SDEs
\begin{equation}\label{E36}
\begin{aligned}
\d u_1&= \left[ -\lambda_1 u_1- \left(u_2^2+ u_3^2\right)\right] \d t +\kappa \d t + \sigma \d W(t),\\
\d u_2&= \left[ -\lambda_2u_2+ u_1u_2\right] \d t, \\
\d u_3&= \left[ -\lambda_3u_3+ u_1u_3\right] \d t,
\end{aligned}
\end{equation}
where $\sigma >0$, $\kappa \in \mathbb{R}$, and $W$ is a standard real-valued  Wiener process.

Clearly we can write the SDEs \REQ{E36} in the form
\begin{equation}\label{E36b}
\d u = \left[Au+B(u,u) + \kappa f_1 \right] \d t + \sigma f_1\d W(t),
\end{equation}
where the maps $A$ and $B$ are defined by
\begin{equation}\label{E36c}
A\left(\begin{array}{c}u_1\\ u_2 \\ u_3\end{array}\right)= -\left(\begin{array}{c}\lambda_1u_1\\ \lambda _2u_2 \\ \lambda_3u_3\end{array}\right),
\quad B\left[ \left(\begin{array}{c}u_1\\ u_2 \\ u_3\end{array}\right), \left(\begin{array}{c} v_1\\ v_2 \\ v_3\end{array} \right)\right]=\left(\begin{array}{c}-u_2v_2-u_3v_3\\ u_2v_1 \\ u_3v_1\end{array}\right)
\end{equation}
and $(f_i)_{i=1}^3$ is the canonical orthonormal basis of $\mathbb{R}^3$.

Note that, as in the case of the NSE-s,  the mapping $B:\mathbb{R}^3\times \mathbb{R}^3 \to \mathbb{R}^3$ is bilinear and
$$
\left(B(u,v),w\right) = b(u,v,w),
$$
where $b$ is a trilinear form on $\mathbb{R}^3$ defined by
$$
b(u,v,w)= -(u_2v_2+ u_3v_3)w_1 + u_2v_1w_2+ u_3v_1w_3.
$$
Note that like for the Navier--Stokes nonlinear mapping  we have
\begin{align}\label{E37}
b(u,v,w)&=-b(u,w,v),\qquad u,v,w\in \mathbb{R}^3,\\
\label{E38}
(B(u,v),v)&= 0,\qquad u,v\in \mathbb{R}^3,\\
\label{E39}
B(f_1,f_1)&=0.
\end{align}
\begin{remark}\label{rem-3.6} Let us emphasize that the condition \eqref{E39} above corresponds to the assumption \eqref{E25} which, as we have pointed out earlier, is satisfied for the 2D Navier--Stokes equations with periodic boundary conditions, see \cite{Flandoli+Gozzi_1998}. Thus our  equation \eqref{E36b} can be seen as  a simple finite dimensional  model of such a problem. Let us point out here that a more general, but still finite dimensional, has been recently investigated by Hairer and Coti-Zelati  in \cite{H+CZ_2020}. They proved the ergodicity of the non-unique invariant probability measures.

 One should also mention a less recent paper \cite{BBOP_2015}, by Ba\v{n}as et al, who studied the uniqueness and non-uniqueness of invariant probability measures for second order stochastic differential equations on a sphere.

 Let us also point out that contrary to 2D  Navier--Stokes equations  we have
\begin{align*}
\label{E39b}
B(f_j,f_j)&=-f_1\not =0  \qquad \text{for  $j=2,3$.}
\end{align*}
\end{remark}

Given $v\in \mathbb{R}^3$ we denote by $u(\cdot;v)$ the solution of \REQ{E36} starting at time $0$ from $v$. Note that
$$
u_1\equiv \kappa/\lambda_1, \quad u_2\equiv 0\equiv u_3,
$$
is a stationary solution to the deterministic problem
\begin{equation}\label{E39b}
\begin{aligned}
\d u_1&= \left[ - \lambda_1u_1- \left(u_2^2+ u_3^2\right) +\kappa  \right]\d t ,\\
\d u_2&= \left[ - \lambda_2u_2+ u_1u_2\right] \d t, \\
\d u_3&= \left[ - \lambda_3u_3+ u_1u_3\right] \d t.
\end{aligned}
\end{equation}
 Note that if $\kappa \le \lambda_1 \min\{\lambda_2, \lambda_3\}$, then there is  unique stationary solution to the system, whereas if $\kappa > \lambda_1 \min\{\lambda_2, \lambda_3\}$, then there exists more than one such a solution. The set of solutions different from the described above can be characterized as follows:

$$
\begin{cases}
\text{If $\lambda _2= \lambda_3$, then $u_1=\lambda_2,\quad u_2^2+u_3^2 = \kappa -\lambda_1\lambda_2$.}\\
\text{If $\lambda_2 >\lambda_3$, $\lambda_2\lambda_3 \ge \kappa$, $\lambda_3\lambda_1<\kappa$ then
$u_1=\lambda_3$, $u_2=0$ and $u^2_3= \kappa -\lambda_1\lambda_3$.}\\
\text{If $\lambda_3 >\lambda_2$, $\lambda_3\lambda_2 \ge \kappa$, $\lambda_2\lambda_1<\kappa$ then
$u_1=\lambda_2$, $u_3=0$ and $u^2_2= \kappa -\lambda_1\lambda_2$.}\\
\text{If $\lambda_1\max\{\lambda_2, \lambda_3\}<\kappa$, $\lambda_2\not = \lambda_3$,
then:}\\
\text{ \qquad (i)$u_1=\lambda_2$, $u_3=0$ and $u^2_2= \kappa -\lambda_1\lambda_2$,}\\
\text{\qquad \quad or}\\
\text{\qquad  (ii)$u_1=\lambda_3$, $u_2=0$ and  $u^2_3= \kappa -\lambda_1\lambda_3$.}
\end{cases}
$$

A natural question arises  whether  the stochastic differential equation \REQ{E36}  exhibits a  similar phenomena as its deterministic counterpart \eqref{E39b}.  We have the following result.
\begin{theorem}\label{T34} In the framework described above the following holds.
\begin{itemize}
\item[$(i)$] For arbitrary parameters, there exists  an invariant probability measure to \REQ{E36}. In fact for any initial value $v\in \mathbb{R}^3$,  there exists  a sequence $t_n\toup +\infty$ such that the following sequence of
Borel probability measures on $\mathbb{R}^3$
$$
 \mathcal{L}\left(\frac{1}{t_n}\int_0^{t_n} u(s;v)\d s\right)
$$
converges weakly to a Borel probability measure on $\mathbb{R}^3$. Consequently, by the Krylov--Bogolyubov theorem,  the simplified  stochastic  NSE \REQ{E36} has at least one invariant probability measure. \\
\item[$(ii)$] For arbitrary $\lambda_1>0$, $\kappa , \sigma\in \mathbb{R}$,  the law $\nu_{\sigma, \kappa}$ of
$$
\frac{ \kappa}{\lambda_1}f_1+\sigma\int_0^{+\infty} \e ^{-\lambda_1t}  \d W(t) f_1
$$
in $(\mathbb{R}^3,\mathcal{B}(\mathbb{R}^3))$ is Gaussian and invariant both for \REQ{E36} and for the stochastic linear "Stokes" equation
\begin{equation}\label{E36d}
dz = Az\d t+ \left(\kappa f_1 \d t + \sigma f_1\d W \right).
\end{equation}
\item[$(iii)$]  If $\kappa< \lambda_1\min\{\lambda_2, \lambda_3\}$, then for any $\sigma\ge 0$, the simplified  stochastic NSE \REQ{E36} admits a unique invariant probability measure $\nu_{\sigma,\kappa}$ that is  stochastically stable;  i.e. for any initial data $v\in \mathbb{R}^3$, the laws $\mathcal{L}(u(t;v))$ converge weakly to $\nu_{\sigma,\kappa}$ as $t\to +\infty$.

\item[$(iv)$] If  $\kappa  >\lambda_1\min\{\lambda_2, \lambda_3\}$,  then  there are invariant probability measures different from the gaussian measure $\nu_{\sigma,\kappa}$.
\end{itemize}
\end{theorem}

\begin{remark}\label{R35}
 We can repeat the first two comments from Remarks \ref{R32}.
 \end{remark}
 \begin{remark}\label{R35a}
We would like to emphasize that the novelty of our results is limited to parts (iii) and (iv) of Theorem \ref{T34}.
\end{remark}
 \begin{remark}\label{R35b}
 We have recently learnt from a talk given by Francesco Morandin about two papers \cite{F+GH+V_2016} and \cite{A+B+C+F+P_2016}, in which infinite dimensional models of NSE-s are studied with the noise acting only on the first mode. Contrary to our case, that model is hypoelliptic and admits a  unique invariant probability measure.
\end{remark}

\delb{
\section{Proof of Theorem \ref{T31}}

\label{sec4}

\begin{proof}[Proof of $(i)$] By the It\^o formula and \REQ{E23} there are $c>0$ and $\rho>0$ such that
\begin{align*}
\mathbb{E}\left\vert u(t;v)\right\vert ^2
&\le \mathbb{E} \left\vert v\right\vert ^2 - \rho \int_0^t  \mathbb{E}\left\Vert u(s;v)\right\Vert ^2 +ct.
\end{align*}
 Thus
$$
\sup_{t>0} \mathbb{E}\, \frac{1}{t}\int_0^t \left\Vert u(s;v)\right\Vert ^2\d s<+\infty.
$$
Consequently, since  the embedding $V\hookrightarrow H$ is compact, the laws of
$$
\frac{1}{t} \int_0^t u(s;v)\d s, \quad t>0
$$
are  tight, and hence, by the Prokhorov theorem, they are   relatively compact in the topology of weak convergence of Borel probability measures on $H$.
\end{proof}
\begin{remark}
Let us point out that the proof of the existence of an  invariant probability measure based on the use of the Krylov--Bogoliubov theorem holds also for the stochastic NSE-s in unbounded
domains, see \cite{BMO_2015} and \cite{Brzezniak+Ferrario_2019}. This is due to an elegant generalisation of the classical Krylov--Bogoliubov theorem in \cite{Maslowski+Seidler_1999}.
\end{remark}

\begin{proof}[Proof of $(ii)$] This part follows directly from \REQ{E24}.
\end{proof}
\begin{proof}[Proof of $(iii)$] Assume now that $v\not \in \text{span}(e)$. Let $z$ be the solution to \REQ{E32} with $z(0)=0$. Consider
$$
y= u(\cdot;v)-z.
$$
Clearly, as $B(z,z)=0$, the process $y$ satisfies
$$
\d y = \left(\mu Ay + B(y+z,y+z)\right)\d t =  \left(\mu Ay + B(z,y) + B(y,z)+B(y,y)\right)\d t,
$$
with initial condition $y(0)=u_0$.
 We have, by \REQ{E23}, that
\begin{align*}
\frac{1}{2}\frac{\d }{\d t} |y(t)|^2 &= -\mu\|y(t)\|^2 + b(y(t),y(t),y(t))+ b(z(t),y(t),y(t)) + b(y(t),z(t),y(t))\\
&= -\mu\|y(t)\|^2 +  b(y(t),z(t),y(t)).
\end{align*}
By the last estimate in $\REQ{E21}$,
$$
\frac{\d }{\d t} |y(t)|^2 \le -2\mu \|y(t)\|^2 + 2C\|z(t)\|  \|y(t)\| |y(t)|.
$$
Therefore
$$
\frac{\d }{\d t} |y(t)|^2 \le -\mu  \|y(t)\|^2 + \frac{C^2}{\mu}\|z(t)\|^2   |y(t)|^2\le \left( \frac{C^2}{\mu}
 \|z(t)\|^2 -\mu\lambda_0\right)|y(t)|^2.
$$
Thus
$$
|y(t)|^2 \le |v|^2 \exp\left\{ \int_0^t \left( \frac{C^2}{\mu}\|z(s)\|^2 -\mu\lambda_0\right)\d s\right\}.
$$
Indeed, by the ergodicity of the process $z$,  and consequently by the Strong Law of Large Numbers, and finally by equality \eqref{E34}  we infer that
$$
\lim_{t\to +\infty}\frac 1 t \int_0^t \|z(s)\|^2 \d s =  \int_{V}\|z\|^2 \nu_{\lambda,\mu,e,\sigma,\kappa }(\d z)=
\left( \frac{\kappa^2}{\lambda^2\mu^2}+\frac{\sigma^2}{2\lambda \mu}\right)\lambda,\quad \text{a.s.,}
$$
where $\nu_{\lambda,\mu,e,\sigma,\kappa}$ is the  (Gaussian) ergodic,  invariant probability measure to the Langevin equation, see \eqref{E33}.  Therefore, since by  assumptions condition   \REQ{E35} is satisfied,  we deduce that
$$
\exp\left\{ \int_0^t \left( \frac{C^2}{\mu}\|z(s)\|^2 -\mu\lambda_0 \right)\d s\right\}\to 0,\quad \text{a.s.,}
$$
and  the desired conclusion follows.
\end{proof}
}

\section{Proof of Theorem \ref{T34}}
Without loss of generality we can assume in the proof that all processes considered here are continuous. 
\begin{proof}[Proof of $(i)$]  We can repeat the argument from the proof of part (i) of Theorem \ref{T31}. Namely,  by the It\^o formula and \REQ{E38} we have
\begin{align*}
\mathbb{E}\left\vert u(t;v)\right\vert ^2 &= \mathbb{E}\left\vert v\right\vert ^2 + 2\mathbb{E} \int_0^t \left[ \langle u(s;v),Au(s;v)\rangle + \langle u(s;v),\kappa f_1\rangle + \frac{\sigma^2}{2} \right]\d s \\
&\le \mathbb{E} \left\vert v\right\vert ^2 - \rho \int_0^t  \mathbb{E}\left\vert u(s;v)\right\vert ^2 +ct,
\end{align*}
where $\rho=\min\{\lambda_1,\lambda_2,\lambda_3\}$ and $c= c(\rho, \sigma,\kappa)$ is independent of $t$. Here $|\cdot|$ stands for the Euclidean norm in $\mathbb R^3$. Thus
$$
\sup_{t>0} \mathbb{E}\left[\frac{1}{t}\int_0^t \left\vert u(s;v)\right\vert ^2\d s\right]<+\infty.
$$
Consequently, the laws of
$$
\frac{1}{t} \int_0^t u(s;v)\d s, \quad t>0,
$$
are  tight  in $\mathbb  R^3$, and hence relatively weakly compact. Therefore, the existence of an invariant probability measure follows from the Krylov--Bogoliubov theorem.
\end{proof}
\begin{proof}[Proof of $(ii)$] This part follows follows immediately from the fact that $B(f_1,f_1)=0$.
\end{proof}
\begin{proof}[Proof of $(iii)$]  Note that
$$
u_i(t;v)= \exp\left\{\int_0^t u_1(s;v)\d s - \lambda_i t\right\}v_i, \qquad i=2, 3,
$$
and
\begin{equation}\label{E51}
u_1(t;v)=\e^{-\lambda_1 t} v_1 - \int_0 ^t \e^{-\lambda_1 (t-s)} X(s;v) \d s + Z(t),
\end{equation}
where
\begin{equation}\label{E52}
Z(t):= \int_0^t \e^{-\lambda_1(t-s)}\left(\kappa \d s + \sigma \d W(s)\right)
\end{equation}
and
$$
X(t;v):= u_2^2(t;v)+ u_3^2(t;v)\ge 0.
$$
Thus, by \eqref{E36}
$$
u_1(t;v) \le \e^{-\lambda_1 t} v_1  + Z(t),
$$
and consequently,
$$
\begin{aligned}
X(t;v)&= \e^{2 \int_0^t u_1(s;v)\d s} \left( \e^{-2\lambda _2 t}v_2^2+ \e^{-2\lambda _3 t}v_3^2\right)\\
&\le  \e^{2 \int_0^t Z(s)\d s} \left( \e^{-2\lambda _2 t}v_2^2+ \e^{-2\lambda _3 t}v_3^2\right) \exp\left\{\frac{2|v_1|}{\lambda _1}\right\}.
\end{aligned}
$$
Clearly
$$
\e^{-2\lambda _2 t}v_2^2+ \e^{-2\lambda _3 t}v_3^2
\le \e^{-2\lambda t}
 \left( v_2^2+ v_3^2\right),
$$
where $\lambda =\min\{\lambda_2, \lambda_3\}>0$, and therefore
$$
X(t;v)\le \e^{2\int_0^t Z(s)\d s -2\lambda t}\left( v_2^2+ v_3^2\right)\exp\left\{\frac{2|v_1|}{\lambda _1}\right\}.
$$

By the law of large numbers
$$
\frac {1}{t} \int_0^t Z(s) \d s\to \frac{\kappa}{\lambda_1}, \qquad \text{$\mathbb{P}$-a.s. as $t\to+\infty$. }
$$
Thus, as $\kappa <\lambda_1\lambda$  we have
$$
\lim_{t\to +\infty} X(t;v)= 0, \qquad \text{$\mathbb{P}$-a.s.}
$$
From the first equation of \eqref{E36} we conclude
\begin{align*}
u_1(t;v)&= \e^{-\lambda_1 (t-T)} u_1(T;v) -\int_T^t \e ^{-\lambda_1(t-s)}X(s;v)\d s + \int_T^{t} \e^{-\lambda_1(t-s)} (\kappa \d s + \sigma \d W(s))\\
&= R(t,T)+ Z(t),
\end{align*}
where
$$
R(t,T;v):= \e^{-\lambda_1 (t-T)} u_1(T;v) -\int_T^t \e ^{-\lambda_1(t-s)}X(s;v)\d s - \int_0^{T} \e^{-\lambda_1(t-s)} (\kappa \d s + \sigma \d W(s)).
$$
Since $R(t,T;v)\to 0$, $\mathbb P$ a.s., as $t\gg T$ and $t,T\to +\infty$ and $Z(t)$ converges in law to
\begin{equation}\label{E53}
\tilde\nu_{\sigma,\kappa}:= \mathcal{N}\left( \frac{\kappa}{\lambda_1}, \frac{\sigma^2}{2\lambda_1}\right)
\end{equation}
 it follows that $u_1(t;v)$ converges  in law to $\tilde  \nu_{\sigma,\kappa}$, and the desired conclusion follows with
\begin{equation}
\label{nu}
\nu_{\sigma,\kappa}:=
 \tilde  \nu_{\sigma,\kappa}\otimes\delta_0 \otimes\delta_0.
\end{equation}
\end{proof}
\begin{proof}[Proof of $(iv)$] Assume that $\lambda_2=  \min\{\lambda_2,\lambda_3\}$. Let   $u$ be the solution to \REQ{E36} with the initial data $u_1(0)=u_3(0)=0$ and $u_2(0)=1$. Then
$$
u_1(t)= -\int_0^t \e ^{-\lambda_1(t-s)} X(s)\d s + Z(t),
$$
where the process $Z$ is defined in \REQ{E52} and
$$
X(t)= \exp\left\{ 2\int_0^t (u_1(s)-\lambda_2)\d s\right\},\;\; t \geq 0.
$$
Note that under the prescribed initial condition we have
\begin{equation}\label{E54}
X(t)=u_2^2(t).
\end{equation}
Since $X\ge 0$ and $\lambda_1>0$, we have
\begin{align*}
\int_0^t u_1(s)\d s  &=  -\int_0^t \int_0^s  \e ^{-\lambda_1(s-r)} X(r)\d r \d s  + \int_0^t Z(s)\d s \\
&=  -\int_0^t \int_r^t  \e ^{-\lambda_1(s-r)} \d s X(r)\d r   + \int_0^t Z(s)\d s\\
&= \int_0^t\left[  - \frac{1}{\lambda_1} \left(1- \e^{-\lambda_1(t-s)}\right) X(s)   + Z(s)\right] \d s\\
&\ge \int_0^t\left[  - \frac{1}{\lambda_1} X(s)   + Z(s)\right] \d s.
\end{align*}
Therefore, we infer that
$$
X(t)\ge \exp\left\{ 2\int_0^t \left(-\frac{1}{\lambda_1}X(s)+ Z(s)-\lambda_2\right)\d s\right\},\;\; t\geq 0.
$$
Next, let us observe that by the law of large numbers for any $\rho$ such that
$$
0<\rho < \frac{\kappa}{\lambda_1}-\lambda_2,
$$
there exists a random variable $\xi$ such that $\mathbb{P}(\xi >0)=1$ and $\mathbb{P}$-a.s
$$
X(t)\ge \xi \exp\left\{ 2\int_0^t \left(-\frac{1}{\lambda_1} X(s)+ \rho\right)\d s\right\}\quad \text{for all $t>0$.}
$$
Thus
$$
X(t) \exp\left\{ \frac{2}{\lambda_1}\int_0^t X(s)\d s\right\}\ge \xi \e^{2\rho t}.
$$
Equivalently
$$
\frac{\d }{\d t} \exp\left\{ \frac{2}{\lambda_1}\int_0^t X(s)\d s\right\}\ge \frac{2}{\lambda_1} \xi \e^{2\rho t},
$$
and  hence
$$
\exp\left\{ \frac{2}{\lambda_1}\int_0^t X(s)\d s\right\} \ge \frac{\xi}{\rho \lambda_1} \left( \e^{2\rho t}-1\right) +1.
$$
Finally, for $t$ large enough we have
$$
\frac{1}{t}\frac{2}{\lambda_1}\int_0^t X(s)\d s\ge \frac{1}{t} \log\left\{ \frac{\xi}{\rho\lambda_1} \left( \e^{2\rho t}-1\right) +1\right\}.
$$
Since
$$
\lim_{t\to +\infty}  \frac{1}{t}  \log\left\{ \frac{\xi}{\rho\lambda_1} \left( \e^{2\rho t}-1\right) +1\right\}= 2\rho
$$
we can  see that
\begin{equation}\label{E55}
\liminf_{t\to +\infty} \frac{1}{t} \frac{2}{\lambda_1}\int_0^t X(s)\d s \ge 2\rho.
\end{equation}
This implies that there exists an invariant probability measure different from $\nu_{\sigma,\kappa}$ defined in \eqref{nu}.  Indeed, consider the Markov process $(u_1,u_2, u_3,  X=u_2^2)$, see \REQ{E54},  for initial value $(0,1,0,1)$. From the first part of the theorem, the sequence of laws
$$
\mathcal{L}\left(\frac 1 t \int_0^t(u_1(s),u_2(s), u_3(t), X(s))\d s\right)
$$
is tight and hence there is a sequence $t_n\toup +\infty$  and a probability  measure $\nu$ on $\mathbb{R}^3\times [0,+\infty)$ such that
$$
\mathcal{L}\left(\frac 1 {t_n} \int_0^{t_n}(u_1(s),u_2(s), u_3(t), X(s))\d s\right)
$$
converge to $\nu$. The probability measure
$$
\tilde \nu(\Gamma)=\nu(\Gamma\times [0,+\infty)), \quad \Gamma \in
\mathcal{B}(\mathbb{R}^3)
$$
is invariant  for  the process $(u_1(t),u_2(t), u_3(t))$, $t\ge 0$. Since, thanks to \REQ{E55}, its marginal with respect to the second variable is not $\delta_0$, it is different from $\nu_{\sigma,\kappa}$.
\end{proof}

\begin{remark}\label{R51}
Let $\tilde\nu_{\sigma, \kappa}$ be given by \REQ{E53}. The method of the proof of part (iii)  of Theorem \ref{T31}  yields  the following criterion: if
\begin{equation}\label{E56}
2\int_{\mathbb{R}}|z| \tilde\nu _{\sigma, \kappa}(\d z)<\min\{\lambda_1,\lambda_2,\lambda_3\},
\end{equation}
then $\nu_{\sigma, \kappa}$, see \eqref{nu}, is the unique invariant probability measure for the nonlinear equation.

This condition is stronger than the condition $\kappa < \lambda_1 \min\{\lambda_2, \lambda_3\}$ for an arbitrary $\sigma$, appearing in Theorem \ref{T31}(iii). In particular, for a fixed $\kappa \ge 0$, \REQ{E56} is violated for large $\sigma$ (see also Remark \ref{R32}).
\end{remark}
To see that \REQ{E56} is really a sufficient condition for ergodicity denote by $z$ the solution of  the linear equation
$$
\d z= Az\d t + \left(\kappa f_1\d t+ \sigma f_1\d W(t)\right), \qquad z(0)=0.
$$
Let $v\in \mathbb{R}^3$. Then $y=u(\cdot;v)- z$ satisfies
$$
\d y = \left[ Ay + B(y+z,y+z)\right]\d t, \qquad y(0)=v.
$$
Hence
\begin{align*}
\frac{1}{2} \frac{\d }{\d t} \left \vert y(t)\right\vert ^2 &= \langle Ay,y\rangle + b(y(t),z(t),y(t)).
\end{align*}
Clearly
$$
 \langle Ay,y\rangle \le -\overline {\lambda}\left\vert y(t)\right\vert ^2,
$$
where
$$
\overline{\lambda}:= \min\{\lambda_1,\lambda_2,\lambda_3\}.
$$
Next, it is easy to see that
$$
\left\vert b(y(t),z(t),y(t))\right\vert \le 2\left\vert z(t)\right\vert \left\vert y(t)\right\vert ^2.
$$
Consequently we have the estimate
$$
\frac{1}{2} \frac{\d }{\d t} \left \vert y(t)\right\vert ^2 \le \left(-\overline {\lambda}+ 2 \left\vert z(t)\right\vert\right)\left\vert y(t)\right\vert ^2,
$$
and hence
$$
\left \vert y(t)\right\vert ^2 \le \left \vert v\right\vert ^2 \exp\left\{ 2\int_0^t \left(-\overline {\lambda}+ 2 \left\vert z(s)\right\vert\right)\d s\right\}.
$$
Since, by the ergodicity of $\tilde \nu_{\sigma, \kappa}$ for $z$,
$$
\frac{1}{t}\int_0^t \left\vert z(s)\right\vert\ \d s\to \int_{\mathbb{R}}|z|\tilde \nu_{\sigma, \kappa}(\d z),
$$
and  the desired conclusion follows.

\medskip
\noindent
\textbf{ Acknowledgements}. 
This work was partially supported by the grant 346300 for IMPAN from the Simons Foundation and the matching 2015-2019 Polish MNiSW fund. T.K.  acknowledges the support of the National Science Centre: NCN grant 2016/23/B/ST1/00492. S.P.  acknowledges the support of the   National Science Center grant 2017/25/B/ST1/02584.  The authors would like to thank the Banach centre for its hospitality.

\end{document}